\documentclass[reqno,centertags, 12pt]{amsart}
\usepackage{amsmath,amsthm,amscd,amssymb}
\usepackage{latexsym}
\newcommand{\bbR}{{\mathbb{R}}}

\newcommand{\bbZ}{{\mathbb{Z}}}
\newcommand{\bbC}{{\mathbb{C}}}

\newcommand{\dott}{\,\cdot\,}

\newcommand{\lb}{\label}
\newcommand{\f}{\frac}

\newcommand{\ti}{\tilde  }

\newcommand{\spec}{\text{\rm{spec}}}

\newcommand{\ess}{\text{\rm{ess}}}

\newcommand{\s}{\text{\rm{s}}}

\newcommand{\supp}{\text{\rm{supp}}}

\newcommand{\bi}{\bibitem}

\newcommand{\beq}{\begin{equation}}
\newcommand{\eeq}{\end{equation}}
\newcommand{\ba}{\begin{align}}
\newcommand{\ea}{\end{align}}
\newcommand{\veps}{\varepsilon}

\newcounter{smalllist}
\newenvironment{SL}{\begin{list}{{\rm\roman{smalllist})}}{%
\setlength{\topsep}{0mm}\setlength{\parsep}{0mm}\setlength{\itemsep}{0mm}%
\setlength{\labelwidth}{2em}\setlength{\leftmargin}{2em}\usecounter{smalllist}%
}}{\end{list}}

\DeclareMathOperator{\sgn}{sgn}
\DeclareMathOperator{\Real}{Re}
\DeclareMathOperator{\ran}{Ran}
\DeclareMathOperator{\Ima}{Im}
\DeclareMathOperator{\Arg}{Arg}

\allowdisplaybreaks
\numberwithin{equation}{section}

\newtheorem{theorem}{Theorem}[section]
\newtheorem*{RT}{Rakhmanov's Theorem}
\newtheorem*{ERT}{Extended Rakhmanov Theorem} 
\newtheorem*{DRT}{Denisov-Rakhmanov Theorem} 

\newtheorem{proposition}[theorem]{Proposition}

\newtheorem{corollary}[theorem]{Corollary}
\theoremstyle{definition}

\theoremstyle{remark}

\newcommand{\abs}[1]{\lvert#1\rvert}

\begin{document}
\title{Sturm Oscillation and Comparison Theorems}
\author{Barry Simon}

\thanks{$^1$ Mathematics 253-37, California Institute of Technology, Pasadena, CA 91125. 
E-mail: bsimon@caltech.edu. Supported in part by NSF grant DMS-0140592} 

\date{October 15, 2003}

\begin{abstract} This is a celebratory and pedagogical discussion of Sturm oscillation 
theory. Included is the discussion of the difference equation case via determinants 
and a renormalized oscillation theorem of Gesztesy, Teschl, and the author. 
\end{abstract}

\maketitle

\section{Introduction} \lb{s1} 

Sturm's greatest contribution is undoubtedly the introduction and focus on 
Sturm-Liouville operators. But his mathematically deepest results are clearly 
the oscillation and comparison theorems. In \cite{Sturm1,Sturm2}, he discussed 
these results for Sturm-Liouville operators. There has been speculation that 
in his unpublished papers he had the result also for difference equations, since 
shortly before his work on Sturm-Liouville operators, he was writing about zeros 
of polynomials, and there is a brief note referring to a never published manuscript 
that suggests he had a result for difference equations. Indeed, the Sturm oscillation 
theorems for difference equations written in terms of orthogonal polynomials are 
clearly related to Descartes' theorem on zeros and sign changes of coefficients. 

In any event, the oscillation theorems for difference equations seem to have 
appeared in print only in 1898 \cite{Bo}, and the usual proof given these days is by 
linear interpolation and reduction to the ODE result. One of our purposes here 
is to make propaganda for the approach via determinants and orthogonal polynomials 
(see Section~\ref{s2}). Our discussion in Section~\ref{s3} and \ref{s4} is more 
standard ODE theory \cite{CL} --- put here to have a brief pedagogical discussion 
in one place. Section~\ref{s5} makes propaganda for what I regard as some interesting 
ideas of Gesztesy, Teschl, and me \cite{GST}. Section~\ref{s6} has three applications 
to illustrate the scope of applicability. 

Our purpose here is celebratory and pedagogical, so we make simplifying assumptions, 
such as only discussing bounded and continuous perturbations. Standard modern 
techniques allow one to discuss much more general perturbations, but this is not 
the place to make that precise. And we look at Schr\"odinger operators, rather 
than the more general Sturm-Liouville operators. 

We study the ODE 
\begin{equation} \lb{1.1} 
Hu=-\f{d^2 u}{dx^2} + Vu =Eu
\end{equation}
typically on $[0,a]$ with $u(0)=u(a)=0$ boundary conditions or on $[0,\infty)$ with 
$u(0)=0$ boundary conditions. The discrete analog is 
\begin{equation} \lb{1.2} 
(hu)_n = a_n u_{n+1} + b_n u_n + a_{n-1} u_{n-1} =Eu
\end{equation}
for $n=1,2,\dots$ with $u_0\equiv 0$. 

\medskip
It is a pleasure to thank W.~Amrein for the invitation to give this talk and for 
organizing an interesting conference, Y.~Last and G.~Kilai for the hospitality 
of Hebrew University where this paper was written, and F.~Gesztesy for useful 
comments. 

\bigskip
\section[Determinants, OP's, and Sturm Theory for Difference Equations]
{Determinants, Orthogonal Polynomials, and Sturm Theory for Difference Equations} 
\lb{s2} 

Given a sequence of parameters $a_1, a_2, \dots$ and $b_1, b_2$ for the difference 
equation \eqref{1.2}, we look at the fundamental solution, $u_n(E)$, defined 
recursively by $u_1 (E)=1$ and 
\begin{equation} \lb{2.1} 
a_n u_{n+1}(E) + (b_n-E) u_n(E) + a_{n-1} u_{n-1}(E)=0 
\end{equation}
with $u_0\equiv 0$, so 
\begin{equation} \lb{2.2} 
u_{n+1}(E) = a_n^{-1} (E-b_n) u_n (E) - a_n^{-1} a_{n-1} u_{n-1}(E) 
\end{equation}
Clearly, \eqref{2.2} implies, by induction, that $u_{n+1}$ is a polynomial of 
degree $n$ with leading term $(a_n \dots a_1)^{-1} E^n$. Thus, we define for 
$n=0,1,2,\dots$ 
\begin{equation} \lb{2.3} 
p_n(E) = u_{n+1}(E) \qquad P_n(E) = (a_1 \dots a_n) p_n(E) 
\end{equation}
Then \eqref{2.1} becomes 
\begin{equation} \lb{2.3a} 
a_{n+1} p_{n+1}(E) + (b_{n+1}-E) p_n(E) + a_n p_{n-1}(E) =0 
\end{equation} 
for $n=0,1,2, \dots$. One also sees that 
\begin{equation} \lb{2.3b} 
EP_n(E) = P_{n+1}(E) + b_{n+1}(E) P_n(E) + a_n^2 P_{n-1}(E) 
\end{equation}
We will eventually see $p_n$ are orthonormal polynomials for a suitable 
measure on $\bbR$ and the $P_n$ are what are known as monic orthogonal 
polynomials. 

Let $J_n$ be the finite $n\times n$ matrix
\[
J_n= \begin{pmatrix}
b_1 & a_1 & 0  \\
a_1 & b_2 & a_2 \\
0 & a_2 & b_3 & \ddots \\
{} & {} & \ddots & \ddots & \ddots\\
{} & {} & {} & \ddots & b_{n-1} & a_{n-1} \\
{} & {} & {} & {} & a_{n-1} & b_n
\end{pmatrix}
\]

\begin{proposition}\lb{P2.1} The eigenvalues of $J_n$ are precisely the 
zeros of $p_n(E)$. We have
\begin{equation} \lb{2.4} 
P_n(E) =\det (E-J_n) 
\end{equation}
\end{proposition} 

\begin{proof} Let $\varphi(E)$ be the vector $\varphi_j(E) = p_{j-1}(E)$,   
$j=1,\dots, n$. Then \eqref{2.1} implies 
\begin{equation} \lb{2.5} 
(J_n-E) \varphi(E) = -a_n p_n (E) \delta_n 
\end{equation}
where $\delta_n$ is the vector $(0,0, \dots, 0,1)$. Thus every zero of $p_n$ is 
an eigenvalue of $J_n$. Conversely, if $\ti\varphi$ is an eigenvector of $J_n$, 
then both $\ti\varphi_j$ and $\varphi_j$ solve \eqref{2.2}, so $\ti\varphi_j 
=\ti\varphi_1 \varphi_j(E)$. This implies that $E$ is an eigenvalue only if $p_n(E)$ 
is zero and that eigenvalues are simple. 

Since $J_n$ is real symmetric and eigenvalues are simple, $p_n(E)$ has $n$ 
distinct eigenvalues $E_j^{(n)}$, $j=1,\dots,n$ with $E_{j-1}^{(n)} < E_j^{(n)}$. 
Thus, since $p_n$ and $P_n$ have the same zeros, 
\[
P_n(E) = \prod_{j=1}^n (E-E_j^{(n)}) = \det (E-J_n)
\]
\end{proof} 

\begin{proposition} \lb{P2.2} 
\begin{SL} 
\item[{\rm{(i)}}] The eigenvalues of $J_n$ and $J_{n+1}$ strictly interlace, that is, 
\begin{equation} \lb{2.6} 
E_1^{(n+1)} < E_1^{(n)} < E_2^{(n+1)} < \cdots < E_n^{(n)} < E_{n+1}^{(n+1)}  
\end{equation}
\item[{\rm{(ii)}}] The zeros of $p_n(E)$ are simple, all real, and strictly interlace 
those of $p_{n+1}(E)$. 
\end{SL} 
\end{proposition} 

\begin{proof} (i) $J_n$ is obtained from $J_{n+1}$ by restricting the quadratic form 
$u\to \langle u, J_{n+1}u\rangle$ to $\bbC^n$, a subspace. It follows that $E_1^{(n+1)} 
=\min_{u,\|u\|=1} \langle u, J_{n+1}u\rangle \leq \min_{u\in\bbC^n, \|u\|=1} \langle 
u,J_{n+1}u\rangle = E_1^{(n)}$. More generally, using that min-max principle 
\[
E_j^{(n+1)} = \max_{\varphi_1, \dots, \varphi_{j-1}} \, \min_{\substack{ \|u\|=1 \\ 
u\perp \varphi_1, \dots, \varphi_{j-1}}} \, \langle u,J_{n+1} u\rangle 
\]
one sees that 
\[
E_j^{(n)} \geq E_j^{(n+1)} 
\]
By replacing $\min$'s with $\max$'s, 
\[
E_j^{(n)} \leq E_{j+1}^{(n+1)}
\]

All that remains is to show that equality is impossible. If $E_0\equiv E_j^{(n)} = E_j^{(n+1)}$ 
or $E_0\equiv E_j^{(n)} = E_j^{(n+1)}$, then $p_{n+1}(E_0) =p_n (E_0)=0$. By \eqref{2.3a}, 
this implies $p_{n-1}(E_0)=0$ so, by induction, $p_0(E)=0$. But $p_0\equiv 1$. Thus 
equality is impossible. 

\smallskip
(ii) Given \eqref{2.4}, a restatement of what we have proven about the eigenvalues of 
$J_n$. 
\end{proof} 

Here is our first version of Sturm oscillation theorems: 

\begin{theorem}\lb{T2.3} Suppose $E_0$ is not an eigenvalue of $J_k$ for $k=1,2,\dots, n$. 
Then 
\begin{align} 
\#(j\mid E_j^{(n)} > E_0) &= \#\{\ell=1,\dots, n\mid \sgn (P_{\ell-1}(E_0)) \neq 
\sgn(P_\ell (E_0))\} \lb{2.7} \\ 
\#(j\mid E_j^{(n)} < E_0) &= \#\{\ell=1, \dots, n\mid \sgn(P_{\ell-1}(E_0) = 
\sgn (P_\ell (E_0)) \} \lb{2.8}
\end{align} 
\end{theorem} 

\begin{proof} \eqref{2.7} clearly implies \eqref{2.8} since the sum of both sides of the 
equalities is $n$. Thus we need only prove \eqref{2.7}. 

Suppose that $E_1^{(\ell)} <\cdots < E_k^{(\ell)} < E_0 < E_{k+1}^{(\ell)} < E_n^{(\ell)}$. 
By eigenvalue interlacing, $J_{\ell+1}$ has $k$ eigenvalues in $(-\infty, 
E_k^{(\ell)})$ and $n-k$ eigenvalues in $(E_{k+1}^{(\ell)},\infty)$. The question is 
whether the eigenvalue in $(E_k^{(\ell)}, E_{k+1}^{(\ell)})$ lies above $E_0$ or below. 
Since $\sgn\,\det (E-J^{(\ell+1)})=(-1)^{\#(j\mid E_j^{(\ell)}>E_0)}$, and similarly for 
$J_{\ell+1}$, and there is at most one extra eigenvalue above $E_0$, we see  
\begin{align*} 
\sgn P_\ell (E_0) &= \sgn P_{\ell+1}(E_0) \Leftrightarrow \#(j\mid E_j^{(\ell)} > E_0) = 
\#(j\mid E_j^{(\ell+1)} > E_0 )\\
\sgn P_\ell (E_0) &= \sgn P_{\ell+1} (E_0) \Leftrightarrow \#j(\mid E_j^{(\ell)} >E_0) 
+1 = \#(j\mid E_j^{(\ell+1)} > E_0 )
\end{align*}

\eqref{2.7} follows from this by induction. 
\end{proof} 

We want to extend this in two ways. First, we can allow $P_k(z_0)=0$ for some $k<n$. 
In that case, by eigenvalue interlacing, it is easy to see $J_{k+1}$ has one more 
eigenvalue than $J_{k-1}$ in $(E_0, \infty)$ and also in $(-\infty, E_0)$, so $\sgn 
(P_{k-1}(z_0)) =-\sgn (P_{k+1}(z_0))$ (also evident from \eqref{2.3b} and $P_k(z_0)=0$). 
Thus we need to be sure to count the change of sign from $< 0,0$ to $>0,a$ as only a 
simple change of sign. We therefore have 

\begin{proposition}\lb{P2.4} \eqref{2.7} and \eqref{2.8} remain true so long as $P_n(E_0) 
\neq 0$ so long as we define $\sgn(0)=1$. If $P_n(E_0)=0$, they remain true so long 
as $\ell =n$ is dropped from the right side. 
\end{proposition}
 
One can summarize this result as follows: For $x\in [0,n]$, define $y(x)$ by linear 
interpolation, that is, 
\[
x=[x]+(x)\Rightarrow y(x) =P_{[x]} + (x) (P_{[x]+1} - P_{[x]}) 
\]
Then the number of eigenvalues of $J_n$ above $E$ is the number of zeros of 
$y(x,E)$ in $[0,n)$. If we do the same for $\ti y$ with $P_{[x]}$ replaced by 
$(-1)^{[x]} P_{[x]}$, then the number of eigenvalues below $E$ is the number of 
zeros of $\ti y$ in $[0,n)$. Some proofs (see \cite{DS}) of oscillation theory for 
difference equations use $y$ and mimic the continuum proof of the next section. 

The second extension involves infinite Jacobi matrices. In discussing eigenvalues 
of an infinite $J$, domain issues arise if $J$ is not bounded (if the moment 
problem is not determinate, these are complicated issues; see Simon \cite{S270}). 
Thus, let us suppose 
\begin{equation} \lb{2.9} 
\sup_n \, (\abs{a_n} + \abs{b_n}) <\infty 
\end{equation}

If $J$ is bounded, the quadratic form of $J_n$ is a restriction of $J$ to $\bbC^n$. 
As in the argument about eigenvalues interlacing, one shows that if $J$ has only 
$N_0 <\infty$ eigenvalues in $(E_0, \infty)$, then $J_n$ has at most $N_0$ 
eigenvalues there. Put differently, if $E_1^{(\infty)} >E_2^{(\infty)} >\cdots$ 
are the eigenvalues of $J$, $E_j^{(\infty)}\geq E_j^{(n)}$. Thus, if $N_n(E) =\#$ 
of eigenvalues of $J_n$ in $(E,\infty)$ and $N_\infty$ the dimension of $\ran 
P_{(E,\infty)}(J)$, the spectral projection 
\begin{equation} \lb{2.10} 
N_n(E) \leq N_{n+1}(E)\leq \dots \leq N_\infty(E) 
\end{equation}

On the other hand, suppose we can find an orthonormal set $\{\varphi_j\}_{j=1}^N$ 
with $M_{jk}^{(\infty)} =\langle\varphi_j, J\varphi_k\rangle =e_j \delta_{jk}$ and 
$\min (e_j)=e_0 >E_0$. If $M_{jk}^{(n)}=\langle \varphi_j, J_n \varphi_k\rangle$, 
$M^{(n)}\to M^{(\infty)}$, so for $n$ large, $M^{(n)} \geq \min (e_j) + \f12 
(e_0 - E_0) > E_0$. Thus $N_n (E_0)\geq N$ for $n$ large. It follows that $\lim 
N_n \geq N_\infty$, that is, we have shown that $N_\infty (E_0) =
\lim_{n\to\infty} N_n(E_0)$. Thus, 

\begin{theorem}\lb{T2.5} Let $J$ be an infinite Jacobi matrix with \eqref{2.9}. 
Then {\rm{(}}with $\sgn(0)=1${\rm{)}} we have 
\begin{align} 
N_\infty (E_0) &= \#\{\ell=1,2, \ldots \mid \sgn (P_{\ell-1}(E_0)) \neq \sgn 
(P_\ell (E_0))\} \lb{2.11}  \\
\dim P_{(-\infty, E_0)}(J) &= \#\{\ell=1,2,\ldots \mid \sgn (P_{\ell-1}(E_0)) 
=\sgn (P_\ell (E_0))\} \lb{2.12}
\end{align}
\end{theorem} 

\begin{corollary} \lb{C2.6} $a_-\leq J \leq a_+$ if and only if for all $\ell$, 
\begin{equation} \lb{2.13} 
P_\ell (a_+) >0 \qquad\text{and}\qquad (-1)^\ell P_\ell (a_-) >0 
\end{equation}
\end{corollary} 

While on the subject of determinants and Jacobi matrices, I would be remiss if I 
did not make two further remarks. 

Given \eqref{2.4}, \eqref{2.3b} is an interesting relation among determinants, and 
you should not be surprised it has a determinantal proof. The matrix $J_{n+1}$ has 
$b_{n+1}$ and $a_n$ in its bottom row. The minor of $E-b_{n+1}$ in $E-J_{n+1}$ is 
clearly $\det (E-J_n)$. A little thought shows the minor of $-a_n$ is $-a_n \det 
(E-J_{n-1})$. Thus 
\begin{equation} \lb{2.14} 
\det (E-J_{n+1}) = (E-b_{n+1}) \det (E-J_n) -a_n^2 \det (E-J_{n-1}) 
\end{equation}
which is just \eqref{2.3b}. 

Secondly, one can look at determinants where we peel off the top and left rather 
than the right and bottom. Let $J^{(1)}, J^{(2)}$ be the Jacobi matrices obtained 
from $J$ by removing the first row and column, the first two, $\dots$. Making the 
$J$-dependence of $P_n (\dott)$ explicit, Cramer's rule implies 
\begin{equation} \lb{2.15} 
(z-J_n)_{11}^{-1} = \f{P_{n-1}(z, J^{(1)})}{P_n(z,J)}
\end{equation} 
In the OP literature, $a_1^{-1} p_n (z,J^{(1)})$ are called the second kind polynomials. 

The analog of \eqref{2.14} is 
\[
P_n (z,J) = (z-b_1) P_{n-1}(z,J^{(1)}) - a_1^2 P_{n-2} (z,J^{(2)})
\]
which, by \eqref{2.15}, becomes 
\begin{equation} \lb{2.16} 
[(z-J)_{11}^{-1}]^{-1} = \f{1}{(z-b_1) - a_1^2 (z-J_{n-1}^{(1)})_{11}^{-1}}
\end{equation}
In particular, since $d\gamma$ is the spectral measure of $\delta_1,J$, we have 
\begin{equation} \lb{2.17} 
(z-J)_{11}^{-1} = \int \f{d\gamma(x)}{z-x} \equiv -m (z,J)  
\end{equation}
and \eqref{2.16} becomes in the limit with $(z-J^{(1)})_{11}^{-1} \to -m (z,J^{(1)})$ 
\begin{equation} \lb{2.18} 
m(z;J) = \f{1}{b_1 - z - a_1^2 m(z;J^{(1)})}
\end{equation}
\eqref{2.16} leads to a finite continued fraction expansion of $(z-J_n)_{11}^{-1}$ 
due to Jacobi, and \eqref{2.18} to the Stieltjes continued fraction. Sturm's 
celebrated paper on zeros of polynomials is essentially also a continued fraction 
expansion. It would be interesting to know how much Sturm and Jacobi knew of 
each other's work. Jacobi visited Paris in 1829 (see James \cite{James}), but 
I have no idea if he and Sturm met at that time. 

\bigskip
\section{Sturm Theory of the Real Line} \lb{s3} 

We will suppose $V$ is a bounded function $[0,\infty)$. We are interested in 
solutions of 
\begin{equation} \lb{3.1} 
-u'' + Vu =Eu 
\end{equation}
for $E$ real. 

\begin{theorem}[Sturm Comparison Theorem]\lb{T3.1} For $j=1,2$, let $u_j$ be not 
identically zero and solve $-u''_j + Vu_j = E_j u_j$. Suppose $a<b$, $u_1 (a) = u_1(b) =0$ 
and $E_2 > E_1$. Then $u_2$ has a zero in $(a,b)$. If $E_2 =E_1$ and $u_2 (a) \neq 0$, 
then $u_2$ has a zero in $(a,b)$. 
\end{theorem} 

\begin{proof} Define the Wronskian 
\begin{equation} \lb{3.2} 
W(x) = u'_1 (x) u_2(x) - u_1(x) u'_2 (x)
\end{equation}
Then 
\begin{equation} \lb{3.3} 
W'(x) =(E_2 - E_1) u_1(x) u_2(x) 
\end{equation}

Without loss, suppose $a$ and $b$ are successive zeros of $u_1$. By changing signs 
of $u$ if need be, we can suppose $u_1 >0$ on $(a,b)$ and $u_2>0$ on $(a,a+\veps)$ 
for some $\veps$. Thus $W(a) =u'_1(a) u_2(a)\geq 0$ (and, in case $E_1=E_2$ and $u_2(a) 
\neq 0$, $W(a)>0$). If $u_2$ is nonvanishing in $(a,b)$, then $u_2 \geq 0$ there, so 
$W(b)>0$ (if $E_2 > E_1$, $(E_2 -E_1)\int_a^b u_1 u_2 \, dx >0$, and if $E_2=E_1$ 
but $u_2(a)\neq 0$, $W(a)>0$). Since $W(b)=u'_1(b) u_2(b)$ with $u'_1(b)<0$ and 
$u_2(b)\geq 0$, this is impossible. Thus we have the result by contradiction. 
\end{proof} 

\begin{corollary}\lb{C3.2} Let $u(x,E)$ be the solution of \eqref{3.1} with 
$u(0,E)=0$, $u'(0,E)=1$. Let $N(a,E)$ be the number of zeros of $u(x,E)$ in 
$(0,a)$. Then, if $E_2 >E_1$, we have $N(a,E_2)\geq N(a,E_1)$ for all $a$. 
\end{corollary} 

\begin{proof} If $n=N(a,E_1)$ and $0<x_1 < \cdots < x_n <a$ are the zeros of 
$u(x,E_1)$, then, by the theorem, $u(x,E_2)$ has zeros in $(0,x_1), (x_1, x_2), 
\dots, (x_{n-1}, x_n)$. 
\end{proof} 

This gives us the first version of the Sturm oscillation theorem: 

\begin{theorem}\lb{T3.3} Let $E_0 < E_1 < \cdots$ be the eigenvalues of $H\equiv 
-\f{d^2}{dx^2} +V(x)$ on $L^2 (0,a)$ with boundary conditions $u(0)=u(a)=0$. Then 
$u(x,E_n)$ has exactly $n$ zeros in $(0,a)$. 
\end{theorem} 

\begin{proof} If $u_k\equiv u(\dott, E_k)$ has $m$ zeros $x_1 < x_2 < \cdots x_m$ 
in $(0,a)$, then for any $E>E_k$, $u(\dott, E)$ has zeros in $(0,x), \dots, (x_{m-1}, 
x_m),(x_m, a)$ and so, $u_{k+1}$ has at least $m+1$ zeros. It follows by induction 
that $u_n$ has at least $n$ zeros, that is, $m\geq n$. 

Suppose $u_n$ has $m$ zeros $x_1 <\cdots < x_m$ in $(0,a)$. Let $v_0, \dots, v_m$ 
be the function $u_n$ restricted successively to $(0,x_1), (x_1, x_2), \dots, 
(x_m, a)$. The $v$'s are continuous and piecewise $C^1$ with $v_\ell (0) = v_\ell 
(a)=0$. Thus they lie in the quadratic form domain of $H$ (see \cite{RS1,RS2} for 
discussions of quadratic forms) and 
\begin{align} 
\langle v_j,Hv_k\rangle &= \int_0^a v'_j v'_k + \int_0^a V v_j v_k \notag \\
&=\delta_{jk} E \int_0^a v_j^2\, dx \lb{3.4} 
\end{align} 
since if $j=k$, we can integrate by parts and use $-u'' + Vu=Eu$. 

It follows that for any $v$ in the span of $v_j$'s, $\langle v, Hv\rangle =E\|v\|^2$, 
so by the variational principle, $H$ has at least $m+1$ eigenvalues in $(-\infty, E_n)$, 
that is, $n+1 \geq m+1$.  
\end{proof} 

{\it Remark.} The second half of this argument is due to Courant-Hilbert \cite{CH}. 

\smallskip
If we combine this result with Corollary~\ref{C3.2}, we immediately have: 

\begin{theorem}[Sturm Oscillation Theorem]\lb{T3.4} The number of eigenvalues of $H$ 
strictly below $E$ is exactly the number of zeros of $u(x,E)$ in $(0,a)$. 
\end{theorem} 

As in the discrete case, if $H_a$ is $-\f{d^2}{dx^2} + V(x)$ on $[0,a]$ with 
$u(0)=u(a)=0$ boundary conditions and $H_\infty$ is the operator on $L^2 (0,\infty)$ 
with $u(0)=0$ boundary conditions, and if $N_a (E) = \dim P_{(-\infty, E)} (H_a)$, 
then $N_a (E)\to N_\infty (E)$, so 

\begin{theorem}\lb{T3.5} The number of eigenvalues of $H_\infty$ strictly below $E$,  
more generally $\dim P_{(-\infty, E)}(H)$, is exactly the number of zeros of $u(x,E)$ 
in $(0,\infty)$. 
\end{theorem} 

There is another distinct approach, essentially Sturm's approach in \cite{Sturm1}, 
to Sturm theory on the real line that we should mention. Consider zeros of $u(x,E)$, 
that is, solutions of 
\begin{equation} \lb{3.5} 
u(x(E),E) =0 
\end{equation} 
$u$ is a $C^1$ function of $x$ and $E$, and if $u(x_0, E)=0$, then $u'(x_0, E_0)\neq 0$ 
(since $u$ obeys a second-order ODE). Thus, by the implicit function theorem, for $E$ 
near $E_0$, there is a unique solution, $x(E)$, of \eqref{3.4} near $x_0$, and it 
obeys 
\begin{equation} \lb{3.6} 
\left. \f{dx}{dE}\right|_{E_0} =  -\left. \f{\partial u/\partial E}{\partial u/\partial x} 
\right|_{x=x_0, E=E_0} 
\end{equation}

Now, $v\equiv \partial u/\partial E$ obeys the equation 
\begin{equation} \lb{3.7} 
-v'' + Vv = Ev+u 
\end{equation}
by taking the derivative of $-u'' +Vu=Eu$. Multiply \eqref{3.7} by $u$ and integrate by 
parts from $0$ to $x_0$. Since $v(0)=0$, there is no boundary term at $0$, but there is 
at $x_0$, and we find 
\[
v(x_0) u'(x_0) = \int_0^{x_0} \, \abs{u(x)}^2\, dx
\]
Thus \eqref{3.6} becomes 
\begin{equation} \lb{3.8} 
\f{dx_0}{dE} = -\abs{u'(x_0,E)}^{-2} \int_0^{x_0} \, \abs{u(x,E)}^2 \, dx < 0
\end{equation}

Thus, as $E$ increases, zeros of $u$ move towards zero. This immediately implies the 
comparison theorem. Moreover, starting with $u_n$, the $(n+1)$-st eigenfunction at 
energy $E_n$, if it has $m$ zeros in $(0,a)$ as $E$ decreases from $E_n$ to a value, 
$E'$ below $-\|V\|_\infty$ (where $u(x,E')>0$ has no zeros in $(0,\infty)$), the 
$m$ zeros move out continuously, and so $u(a,E)=0$ exactly $m$ times, that is, $m=n$. 
This proves the oscillation theorem. 

\bigskip
\section{Rotation Numbers and Oscillations} \lb{s4} 

Take the solution $u(x,E)$ of the last section and look at the point 
\[
\pi(x,E) = \binom{u'(x,E)}{u(x,E)} 
\]
in $\bbR^2$. $\pi$ is never zero since $u$ and $u'$ have no common zeros. At most 
points in $\bbR^2$, the argument of $\pi$, that is, the angle $\pi$ makes with 
$\binom{1}{0}$, can increase or decrease. $u$ can wander around and around. But not 
at points where $u=0$. If $u'>0$ at such a point, $\pi$ moves from the lower right 
quadrant to the upper right, and similarly, if $u'<0$, it moves from the upper 
left to lower left. Thus, since $\pi$ starts at $\binom{1}{0}$, we see 

\begin{theorem}\lb{T4.1} If $u(x,E)$ has $m$ zeros in $(0,a)$, then $\Arg \pi (a,E)$ 
{\rm{(}}defined by continuity and $\Arg \pi (0,E)=0${\rm{)}} lies in $(m\f{\pi}{2}, 
(m+1)\f{\pi}{2}]$. 
\end{theorem} 

If $u$ and $v$ are two solutions of $-u'' +Vu=Eu$ with $u(0)=0$, $v(0)\neq 0$, we can 
look at 
\[
\ti\pi(x,E) =\binom{u}{v}
\]
$\ti\pi$ is never zero since $u$ and $v$ are linear independent. 
$W(x)=u'v-v'u$ is a constant, say $c$. $c\neq 0$ since $u$ and $v$ are linear 
independent. Suppose $c>0$. Then if $u(x_0)=0$, $u'(x_0) =c/v(x_0)$ has the same 
sign as $v(x_0)$. So the above argument applies (if $c<0$, there is winding in the 
$(u,v)$-plane in the opposite direction). Rather than look at $\ti\pi$, we can look 
at $\varphi=u+iv$. Then $u'v -vu' =\Ima (\bar\varphi \varphi')$. Thus we have 

\begin{theorem}\lb{T4.2} Let $\varphi(x,E)$ obey $-\varphi'' +V\varphi=E\varphi$ and be 
complex-valued with 
\begin{equation} \lb{4.1} 
\Ima (\bar\varphi(0) \varphi'(0))>0 
\end{equation}
Suppose $\Real \varphi(0)=0$. Then, if $\Real\varphi$ has $m$ zeros in $(0,a)$, then 
$\Arg (\varphi(a))$ is in $(m\f{\pi}{2}, (m+1)\f{\pi}{2}]$. 
\end{theorem} 

The ideas of this section are the basis of the relation of rotation numbers and density 
of states used by Johnson-Moser \cite{JM} (see also \cite{Jo}). We will use them as the 
starting point of the next section. 

\bigskip
\section{Renormalized Oscillation Theory} \lb{s5} 

Consider $H=-\f{d^2}{dx^2}+V$ on $[0,\infty)$ with $u(0)=0$ boundary conditions where, as 
usual, for simplicity, we suppose that $V$ is bounded. By Theorem~\ref{T3.5}, $\dim 
P_{(-\infty,E)}(H)$ is the number of zeros of $u(x,E)$ in $(0,\infty)$. If we want to 
know $\dim P_{[E_1, E_2)}(H)$, we can just subtract the number of zeros of $u(x,E_1)$ on 
$(0,\infty)$ from those of $u(x,E_2)$. At least, if $\dim P_{(-\infty, E_2)}(H)$ is 
finite, one can count just by subtracting. But if $\dim P_{(-\infty, E_1)}(H)=\infty$ 
while $\dim P_{[E_1, E_2)}$ is finite, both $u(x,E_2)$ and $u(x,E_1)$ have infinitely many 
zeros, and so subtraction requires regularization. 

One might hope that 
\begin{equation} \lb{5.1} 
\dim P_{[E_1, E_2)}(H) =\lim_{a\to\infty} (N(E_2,a) -N(E_1, a)) 
\end{equation}
where $N(E,a)$ is the number of zeros of $u(x,E)$ in $(0,a)$. This is an approach of 
Hartmann \cite{Hart}. \eqref{5.1} cannot literally be true since $N(E_2, a) -N(E_1,a)$ 
is an integer which clearly keeps changing when one passes through a zero of $u(x,E_2)$ 
that is not also a zero of $u(x,E_1)$. One can show that for $a$ large, the absolute 
value of difference of the two sides of \eqref{5.1} is at most one, but it is not obvious 
when one has reached the asymptotic region. 

Instead, we will describe an approach of Gesztesy, Simon, and Teschl \cite{GST}; see Schmidt 
\cite{Sch} for further discussion. Here it is for the half-line (the theorem is true in much 
greater generality than $V$ bounded and there are whole-line results). 

\begin{theorem}\lb{T5.1} Let $V$ be bounded and let $H=-\f{d^2}{dx^2} +V(x)$ on 
$[0,\infty)$ with $u(0)=0$ boundary conditions. Fix $E_1 <E_2$. Let 
\begin{equation} \lb{5.2} 
W(x) = u(x,E_1) u'(x,E_2) -u' (x,E_1) u(x,E_2) 
\end{equation} 
and let $N$ be the number of zeros of $W$ in $(0,\infty)$. Then 
\begin{equation} \lb{5.3} 
\dim P_{(E_1, E_2)} (H)=N
\end{equation} 
\end{theorem} 

The rest of this section will sketch the proof of this theorem under the assumption that  
$\dim P_{(-\infty, E_2)}(H)=\infty$. This will allow a simplification of the argument 
and covers cases of greatest interest. Following \cite{GST}, we will prove this in three steps: 
\begin{SL} 
\item[(1)] Prove the result in a finite interval $[0,a]$ in case $u(a,E_2)=0$. 
\item[(2)] Prove $\dim P_{(E_1,E_2)}(H)\leq N$ by limits from (1) when $\dim 
P_{(-\infty, E_2)}(H)=\infty$. 
\item[(3)] Prove $\dim P_{(E_1, E_2)}(H)\geq N$ by a variational argument. 
\end{SL} 

\smallskip
\noindent{\bf Step 1.} \ We use the rotation number picture of the last section. Define the 
Pr\"ufer angle $\theta (x,E)$ by 
\begin{equation} \lb{5.4} 
\tan(\theta(x,E)) = \f{u(x,E)}{u'(x,E)} 
\end{equation}
with $\theta(0,E)=0$ and $\theta$ continuous at points, $x_0$, where $u'(x_0, E)=0$. 
Using $\f{d}{dy}\tan y =1+\tan^2 y$, we get 
\begin{equation} \lb{5.5} 
\f{d\theta}{dx} = \f{(u')^2 - uu''}{u^2 + (u')^2} 
\end{equation} 

Let $\theta_1, \theta_2$ be the Pr\"ufer angles for $u_1(x)\equiv u(x,E_1)$ and $u_2 (x) 
\equiv u(x,E_2)$. Suppose $W(x_0)=0$. This happens if and only if $u(x_0,E)/u'(x_0,E_1) 
= u(x_0, E_2)/u'(x_0, E_2)$, that is, $\theta_2=\theta_1 + k\pi$ with $k\in\bbZ$. If it 
happens, we can multiply $u_2$ by a constant so $u_1 (x_0)=u_2 (x_0)$, $u'_1(x_0) = 
u'_2(x_0)$. Once we do that, \eqref{5.5} says 
\[
\f{d}{dx}\, (\theta_2 - \theta_1) = \f{(E_2 - E_1) u_1^2 (x_0)}
{u'_1 (x_0)^2 +u_1^2 (x_0)} > 0
\]
Thus 
\begin{equation} \lb{5.6} 
\theta_1 = \theta_2 \mod \pi \Rightarrow \theta'_2 > \theta'_1 
\end{equation}

Think of $\theta_2$ as a hare and $\theta_1$ as a tortoise running around a track of 
length $\pi$. There are two rules in their race. They can each run in either direction, 
except they can only pass the starting gate going forward (i.e., $\theta_j=0 \mod
\pi\Rightarrow \theta'_j >0$), and the hare can pass the tortoise, not vice-versa 
(i.e., \eqref{5.6} holds). 

Suppose $H_a$, the operator on $(0,a)$ with $u(0)=u(a)=0$ boundary condition, has $m$ 
eigenvalues below $E_2$ and $n$ below $E_1$. Since $u(a,E_2)=0$, $\theta_2(a) = (m+1) 
\pi$, that is, at $x=a$, the hare makes exactly $m+1$ loops of the track. At $x=a$, 
the tortoise has made $n$ loops plus part, perhaps all, of an additional one. 
Since $\theta'_2 -\theta'_1 >0$ at $x=0$, the hare starts out ahead. Thus, the hare 
must overtake the tortoise exactly $m-n$ times between $0$ and $a$ (if $\theta_1(a) 
=(n+1)\pi$, since then $\theta'_2 -\theta'_1 >0$ at $x=0$, $\theta_2 -(m+1)\pi < 
\theta_1 - (n-1)\pi$, and $x=a$; so it is still true that there are exactly $m-n$ 
crossings). Thus 
\begin{equation} \lb{5.7}
P_{(E_1,E_2)}(H_a) =\#\{x_0\in (0,a)\mid W(x_0)=0\} 
\end{equation} 

\smallskip
\noindent{\bf Step 2.} \ Since $\dim P_{(-\infty, E_2)}(H)=\infty$, there is, 
by Theorem~\ref{T3.5}, an infinite sequence $a_1<a_2<\cdots \to \infty$ so that 
$u(a_j, E_2)=0$. $H_{a_j}\to H$ in strong resolvent sense, so by a simple argument, 
\begin{align} 
\dim P_{(E_1, E_2)}(H) &\leq \liminf \dim P_{(E_1, E_2)}(H_a) \notag \\ 
&= N \lb{5.8}
\end{align} 
with $N$ the number of zeros of $W$ in $(0,\infty)$. \eqref{5.8} comes from \eqref{5.7}. 

\smallskip
\noindent{\bf Step 3.} \ Suppose $N<\infty$. Let $0<x_1 < \dots < x_N$ be the zeros of 
$W$\!. Define 
\begin{align} 
\eta_j(x) &= \begin{cases} u_1(x) -\gamma_j u_2(x) & 0<x\leq x_j \\
0 & x\geq x_j \end{cases} \lb{5.9} \\
\notag \\
\ti\eta_j (x) &= \begin{cases} u_1(x) + \gamma_j u_2(x) & 0<x < x_j \\
0 & x > x_j \end{cases}\lb{5.10} 
\end{align} 
where $u_j(x)=u(x,E_j)$ and $\gamma_j$ is chosen by 
\begin{equation} \lb{5.11}
\gamma_j = \begin{cases} u_1(x_j)/u_2(x_j) &\text{if } u(x_j)\neq 0 \\
u'_1 (x_j)/u'_2(x_j) & \text{if } u(x_j)=0 \end{cases}
\end{equation} 
Since $W(x_j)=0$, $\eta_j$ is a $C^1$ function of compact support and piecewise $C^2$, 
and so in $D(H)$. But $\ti\eta$ is discontinuous. 

We claim that if $\eta$ is in the span of $\{\eta_j\}_{j=1}^N$, then 
\begin{equation} \lb{5.12} 
\biggl\|\bigg(H - \f{E_2 + E_1}{2}\biggr)\eta\biggr\| = \f{\abs{E_2 - E_1}}{2}\, 
\|\eta\|
\end{equation} 
Moreover, such $\eta$'s are never a finite linear combination of eigenfunctions of 
$H$. Accepting these two facts, we note that since the $\eta_j$ are obviously linear 
independent, \eqref{5.12} implies $\dim P_{(E_1, E_2)}(H)\geq N$. This, together with 
\eqref{5.8}, proves the result. 

To prove \eqref{5.12}, we note that 
\begin{equation} \lb{5.13}
\biggl( H-\f{E_2 + E_1}{2}\biggr) \eta_j = -\f{\abs{E_2 - E_1}}{2}\, \ti\eta_j
\end{equation}
Since $\ti\eta_j$ is not $C^1$ at $x_j$, no $\ti\eta$ is in $D(H)$, hence no $\eta$ can be 
in $D(H^2)$ (so we get control of $\dim P_{(E_1, E_2)}(H)$, not just $\dim P_{[E_1,E_2]} 
(H)$). 

Next, note that since $W'(x) =(E_2 -E_1) u_2 u_1$, we have if $W(x_i)=W(x_{i+1})=0$ that 
\[
\int_{x_i}^{x_{i+1}} u_1(x) u_2(x)\, dx =0
\]
for $i=0,1,2,\dots, N$ where $x_0=0$. Thus 
\begin{equation} \lb{5.14}
\langle\eta_i, \eta_j\rangle = \langle \ti\eta_i, \ti\eta_j\rangle 
\end{equation}
since if $i<j$, the difference of the two sides is $2(\gamma_i + \gamma_j) \int_{x_i}^{x_j} 
u_1(x) u_2(x)=0$. \eqref{5.14} and \eqref{5.13} implies \eqref{5.12}. That completes the 
proof if $N<\infty$. 

If $N$ is infinite, pick $0<x_1 < \cdots < x_L$ successive zeros and deduce $\dim 
P_{(E_1, E_2)}(H)\geq L$ for all $L$. \hfill\qed 

\bigskip
\section{Some Applications} \lb{s6} 

We will consider three typical applications in this section: one classical (i.e., fifty 
years old!), one recent to difference equations, and one of Theorem~\ref{T5.1}. 

\medskip
\noindent{\bf Application 1: Bargmann's Bound.} \ Let $u$ obey $-u'' +Vu=0$ with 
$u(0)=0$ so, if $V$ is bounded, $u(x)/x$ has a finite limit as $x\downarrow 0$. Also 
suppose $V\leq 0$. 

Define $\ti m=-u'/u$ so 
\begin{equation} \lb{6.1}
\ti m' = \abs{V} + \ti m^2
\end{equation}
since $-V=\abs{V}$. Thus $\ti m$ is monotone increasing. It has a pole at each zero, 
$x_0=0$, $x_1, x_2, \dots, x_\ell, \dots$ of $u$. Define  
\begin{equation} \lb{6.2x} 
b(x) = -\f{xu'}{u} = x\ti m(x) 
\end{equation} 
Then $b(x)$ has limit $-1$ as $x\downarrow 0$ and 
\begin{equation} \lb{6.2}
b'(x) =x\abs{V} + \f{(b+b^2)}{x}
\end{equation}
In particular, 
\begin{equation} \lb{6.3}
-1\leq b \leq 0 \Rightarrow b'(x)\leq x\abs{V}
\end{equation}

By the monotonicity of $\ti m$, there are unique points $0<z_1 < x_1 < \cdots < 
x_{\ell-1} < z_\ell < x_\ell$ where $b_\ell =0$, and since $b\to-\infty$ as $x
\downarrow x_j$, there are last points $y_j \subset [x_{j-1}, z_j]$ where $b(y)=-1$ 
for $j=2,3,\dots, \ell$ and at $y_1 =0$, $b(0)=-1$. Integrating $b'$ from 
$y_j$ to $z_j$, using \eqref{6.3}, we find 
\[
\int_{y_j}^{z_j} x\abs{V(x)}\, dx \geq 1 
\]
so 
\[
\int_0^{x_\ell} x\abs{V(x)}\, dx \geq \ell
\]
By the oscillation theorem, if $N(V)=\dim P_{(-\infty, 0)}(H)$, then 
\begin{equation} \lb{6.4}
N(V) \leq \int_0^\infty x\abs{V(x)}\, dx
\end{equation}
This is Bargmann's bound \cite{Barg}. For further discussion, see Schmidt \cite{Sch2}. 

\smallskip
\noindent{\bf Application 2: Denisov-Rakhmanov Theorem.} \ Rakhmanov \cite{Rakh77,Rakh87} 
(see also \cite{MNT85a}) proved a deep theorem about orthogonal polynomials on the unit circle  
that translates to 

\begin{RT}\lb{RT} If $J$ is an infinite Jacobi matrix , $d\mu =f\, dx + d\mu_\s$ and $f(x) 
>0$ and $x\in [-2,2]$ and $\supp(d\mu_\s)\subset [-2,2]$ {\rm{(}}i.e., $\spec(J)\subset 
[-2,2]${\rm{)}}, then $a_n\to 1$, $b_n\to 0$. 
\end{RT} 

From the 1990's, there was some interest in extending this to the more general result, 
where $\spec (J)\subset [-2,2]$ is replaced by $\ess\,\spec (J)\subset [-2,2]$. By using 
the ideas of the proof of Rakhmanov's theorem, one can prove: 

\begin{ERT} There exist $C(\veps)\to 0$ as $\veps\downarrow 0$ so that if $d\mu =f\, 
dx +du$ and $f(x) >0$ a.e.~$x$ in $[-2,2]$ and $\spec(J) \subset [-2-\veps, 2+C]$, then 
\[
\limsup (\abs{a_n-1} + \abs{b_n})\leq C(\veps)
\]
\end{ERT}

Here is how Denisov \cite{DenPAMS} used this to prove 

\begin{DRT} If $d\mu = f(x)\, dx + d\mu_0$, $f(x)>0$ a.e.~$x\in [-2,2]$ and 
$\sigma_\ess(J)\subset [-2,2]$, then $a_n\to 1$ and $b_n\to 0$. 
\end{DRT} 

His proof goes as follows. Fix $\veps$. Since $J$ has only finitely many eigenvalues in 
$[2+\veps, \infty)$, $P_n (2+\veps)$ has only finitely many sign changes. Similarly, 
$(-1)^n P_n (-2-\veps)$ has only finitely many sign changes. Thus, we can find $N_0$ 
so $P_n (2+\veps)$ and $(-1)^n P_n (-2-\veps)$ both have fixed signs if $n>N_0$. Let 
$\ti a, \ti b$ be given by 
\[
\ti a_n = a_{N_0+n} \qquad \ti b_n =b_{N_0 +n} 
\]
By a use of the comparison and oscillation theorems, $\ti J$ has no eigenvalues in 
$(-\infty, -2-\veps)\cup (2+\veps, \infty)$. Thus, by the Extended Rakhmanov Theorem, 
\[
\limsup (\abs{a_n-1} + \abs{b_n})  = \limsup (\abs{\ti a_n-1} + \abs{\ti b_n}) 
\leq C(\veps)
\]
Since $\veps$ is arbitrary, the theorem is proven. 

\medskip
\noindent{\bf Application 3: Teschl's Proof of the Rofe-Beketov Theorem.} \ Let $V_0(x)$ 
be periodic and continuous. Let $H_0 =-\f{d^2}{dx^2} + V_0$ on $L^2 (0,\infty)$ with 
$u(0)=0$ boundary condition. Then 
\[
\sigma_\ess (H_0) = \bigcup_{j=1}^\infty \, [a_j, b_j]
\]
with $b_j <a_{j+1}$. (In some special cases, there is only a finite union with one 
infinite interval.) $(b_j, a_{j+1})$ are called the gaps. In each gap, $H_0$ has either 
zero or one eigenvalue. Suppose $X(x)\to 0$ as $x\to\infty$, and let $H=H_0 + X$. Since 
$\sigma_\ess (H)=\sigma_\ess (H_0)$, $H$ also has gaps in its spectrum. When is it 
true that each gap has at most finitely many eigenvalues? Teschl \cite{Te1,Te2} 
has proven that if $\int_0^\infty x \abs{X(x)}\, dx <\infty$, then for each $j$, 
the Wronskian, $w(x)$, of $u(x,b_j)$ and $u(x,a_{j+1})$ has only finitely many zeros. 
He does this by showing for $H_0$ that $\abs{X(x)}\to \infty$ as $x\to\infty$ and 
by an ODE perturbation argument, this implies $\abs{w(x)}\to\infty$ for $H$. Thus, by 
the results of Section~\ref{s5}, there are finitely many eigenvalues in each gap.  

It is easy to go from half-line results to whole-line results, so Teschl proves 
if $\int\abs{x}\, \abs{X(x)}\, dx <\infty$, each gap has only finitely many eigenvalues. 

This result was first proven by Rofe-Beketov \cite{RB} with another simple proof in 
Gesztesy-Simon \cite{GS}; see that later paper for additional references. Teschl's 
results are stated for the discrete (Jacobi) case (and may be the first proof for 
the finite difference situation), but his argument translates to the one above for 
Schr\"odinger operators.

\bigskip


\end{document}